\documentclass[11pt,reqno,a4paper]{amsart}

\usepackage[plain]{fullpage}

\usepackage{calrsfs}
\usepackage[OT2,T1]{fontenc}
\usepackage[english]{babel}
\usepackage{newtxtext}
\usepackage{newtxmath}
\usepackage{bbm}
\usepackage{tikz}

\usepackage{etoolbox}
\patchcmd{\section}{\scshape}{\bfseries}{}{}
\makeatletter
\renewcommand{\@secnumfont}{\bfseries}
\makeatother

\theoremstyle{plain}
\newtheorem{theorem}{Theorem}[section]
\newtheorem{lemma}[theorem]{Lemma}
\theoremstyle{definition}
\newtheorem{definition}[theorem]{Definition}
\theoremstyle{remark}
\newtheorem{remark}[theorem]{Remark}

\counterwithin{equation}{section}

\usepackage[unicode=true,pdfusetitle, bookmarks=true,bookmarksnumbered=false,bookmarksopen=false, breaklinks=true,pdfborder={0 0 0},pdfborderstyle={},backref=false,colorlinks=true]{hyperref}
\definecolor{myblue}{rgb}{0.09,0.32,0.44} 
\hypersetup{pdfborder={0 0 0},pdfborderstyle={},colorlinks=true,linkcolor=myblue,citecolor=myblue,urlcolor=blue}
\usepackage{mathtools}
\usepackage{mathrsfs}

\newcommand{\Z}{\mathbf{Z}}
\newcommand{\R}{\mathbf{R}}
\newcommand{\Q}{\mathbf{Q}}
\newcommand{\FF}{\mathbf{F}}
\newcommand{\Fbar}{\overline{\mathbf{F}}}
\newcommand{\Cc}{\mathbf{C}}
\newcommand{\mb}{\mathbbm}
\newcommand{\mc}{\mathcal}
\newcommand{\mr}{\mathrm}
\renewcommand{\geq}{\geqslant}
\renewcommand{\leq}{\leqslant}
\renewcommand{\ge}{\geqslant}
\renewcommand{\le}{\leqslant}

\newcommand{\intersect}{\cap}

\newcommand{\union}{\cup}
\newcommand{\Union}{\bigcup}

\DeclarePairedDelimiter\abs{\lvert}{\rvert}%

\makeatletter
\let\oldabs\abs
\def\abs{\@ifstar{\oldabs}{\oldabs*}}

\begin{document}

\date{\today}
\title{Irreducibility of Littlewood polynomials of special degrees}

\author[L.~Bary-Soroker]{Lior Bary-Soroker}
\address{\normalfont LBS: Raymond and Beverly Sackler School of Mathematical Sciences\\
        	 Tel Aviv University\\
        	 Tel Aviv 69978, Israel}
\email{{\tt barylior@tauexe.tau.ac.il}}

\author[D.~Hokken]{David Hokken}
\address{\normalfont DH: Mathematical Institute\\
        Utrecht University \\
        3508 TA Utrecht, The Netherlands}
\email{{\tt d.p.t.hokken@uu.nl}}

\author[G.~Kozma]{Gady Kozma}
\address{\normalfont GK: Department of Mathematics\\
		The Weizmann Institute of Science\\ 
		Rehovot 76100, Israel}
\email{{\tt gady.kozma@weizmann.ac.il}}

\author[B.~Poonen]{Bjorn Poonen}
\address{\normalfont BP: Department of Mathematics\\
		Massachusetts Institute of Technology\\ 
		Cambridge, MA 02139, USA}
\email{{\tt poonen@math.mit.edu}}


\date{December 26, 2023}

\begin{abstract} 
Let $f$ be sampled uniformly at random from the set of degree~$n$ polynomials whose coefficients lie in $\{ \pm 1\}$. A folklore conjecture, known to hold under GRH, states that the probability that $f$ is irreducible tends to $1$ as $n$ goes to infinity.  We prove unconditionally that
\begin{equation*}
\limsup_{n \to \infty} \mb{P}(f \text{ is irreducible}) = 1.
\end{equation*}
\end{abstract}

\maketitle

\section{Introduction}
Let $f(X) = \sum_{i=0}^n \pm X^i$ be a Littlewood polynomial of degree $n$ sampled uniformly at random; that is, its coefficients are independent random variables taking the values $\pm 1$ with probability $1/2$ each. 
A folklore conjecture~\cite{BKK,BK,BV,Konyagin,OP,Poonen} asserts that
\begin{equation} \label{eq:con}
\lim_{n \to \infty} \mb{P}(f \text{ is irreducible}) = 1.
\end{equation}
Breuillard and Varj\'u \cite{BV} proved that the Generalized Riemann Hypothesis implies \eqref{eq:con}. Together with Koukoulopoulos, the first and third authors \cite[Theorem~3.5]{BKK} proved unconditionally  that 
\begin{equation*}
\liminf_{n \to \infty} \mb{P}(f \textnormal{ is irreducible}) > 0.
\end{equation*}
The goal of this paper is to prove \eqref{eq:con} unconditionally when restricting to a subsequence of degrees:
\begin{theorem} \label{thm:short}
In the notation above,
\begin{equation*}
\limsup_{n \to \infty} \mb{P}(f \textnormal{ is irreducible}) = 1.
\end{equation*}
\end{theorem}
If $n=p-1$ for a prime $p$ such that $2$ generates $(\Z/p\Z)^{\times}$, then $f$ is always irreducible over $\FF_2$, so $\mb{P}(f \textnormal{ is irreducible})=1$ for such $n$. Unfortunately, the infinitude of such primes is the content of the Artin primitive root conjecture, which is open. 
However, similar in spirit to this observation is the following more precise result, which implies Theorem~\ref{thm:short} immediately. 

\begin{theorem} \label{thm:main}
There exists an absolute constant $c>0$ such that the following holds.
Suppose that either $p=2$, or $p\geq 7$ is a prime number such that $2$ generates $(\Z/p^2\Z)^{\times}$.  Let $r\geq 1$ be an integer. Set $n = p^r-1$. Then 
\begin{equation}
\label{E:P(f is irreducible)}
\mb{P}(f \textnormal{ is irreducible}) \geq 1-n^{-c}.
\end{equation}
\end{theorem}

\begin{remark}
We may assume that $f$ is monic, since $f$ is irreducible if and only if $-f$ is.
\end{remark}

\begin{remark}
For every $n \geq 1$, there exists at least one irreducible Littlewood polynomial (see, e.g., \cite[Corollary 2.2]{Martin}), so we may adjust $c$ to make \eqref{E:P(f is irreducible)} hold for any finite list of positive integers $n$.  Thus, in proving Theorem~\ref{thm:main}, we may assume that $n$ is sufficiently large, and it does not matter if $n^{-c}$ is replaced by $O(n^{-c})$.
\end{remark}

\subsection*{Acknowledgements}
We thank Gunther Cornelissen, Arno Fehm, Berend Ringeling, and Mieke Wessel for feedback on an earlier version of this paper. LBS was partially supported by the Israel Science Foundation grant no.~702/19. DH was supported by the Dutch Research Council (NWO), project number OCENW.M20.233. GK was partially supported by the Israel Science Foundation grant no.~607/21. BP was supported in part by National Science Foundation grant DMS-2101040 and Simons Foundation grants \#402472 and \#550033.  DH thanks his hosts at Tel Aviv University for their hospitality and support during his visit in Spring 2023 that led to this project.

\section{Proof of Theorem~\ref{thm:main} for $p=2$}
We assume that $n$ is large and $f$ is monic.
We have $(X-1)f(X) \equiv X^{n+1}-1 \equiv (X-1)^{n+1} \pmod 2$, since $n+1$ is a power of $2$.
Let $g(X) = f(X+1)$, so $g(X) = X^n+\sum_{i=0}^{n-1} g_i X^i$ for some $g_i \in 2\Z$.
The maps $f \mapsto (f \bmod 4) \mapsto (g \bmod 4)$ are injective, so the composition defines a bijection from the set of $2^n$ monic Littlewood polynomials to the set of $2^n$ monic degree~$n$ polynomials in $(\Z/4\Z)[X]$ reducing to $X^n$ in $(\Z/2\Z)[X]$.
Thus the $(g_i \bmod 4)$ take the values $0$ and $2$ uniformly and independently.
Fix $\theta$ as in \cite[Corollary~1]{BKK}, with $\theta \in (0,1/2)$.
With high probability, there exists $i< \theta n$ with $g_i \equiv 2 \pmod 4$; choose the smallest such $i$.
Then the $2$-adic Newton polygon of $g$ has a segment of width $n-i$ and height $1$, so $g$ has a $\Q_2$-irreducible factor of degree $\geq n-i$, and hence a $\Q$-irreducible factor of degree $\geq n-i$, and so does $f$; any other irreducible factor of $f$ has degree $\leq i < \theta n$ (see for example \cite[Section 7.4]{Gou} for a discussion on Newton polygons).
On the other hand, by \cite[Corollary~1(a)]{BKK}, with probability $\geq 1 - n^{-c}$, the polynomial $f$ has no irreducible factors of degree $< \theta n$, so then $f$ is irreducible.

\section{Proof of Theorem~\ref{thm:main} for $p>1470$}
\label{sec:largep}

For any $f \in \Z[X]$ and prime $p$, let $f_p \coloneqq (f \bmod p) \in \FF_p[X]$.
Define a probability measure $\mu$ on $\Z$ by $\mu(1)=\mu(-1)=1/2$; it induces a probability measure on the set of polynomials of degree $n$ with integer coefficients by sampling each coefficient independently according to $\mu$. Write $e(x) \coloneqq e^{2\pi i x}$, and define the Fourier transform $\hat{\mu} \colon \R/\Z \to \Cc$ by
\begin{equation} \label{eq:FT}
\hat{\mu}(\alpha) \coloneqq \sum_{k \in \Z} \mu(k) e(\alpha k) = \frac{e(\alpha) + e(-\alpha)}{2} = \cos(2\pi\alpha).
\end{equation}

\begin{proof}[Proof of Theorem~\ref{thm:main} for $p > 1470$.]
Let $\Phi_m$ denote the $m$th cyclotomic polynomial, viewed in $\FF_2[X]$.
By \cite[Exercise 12, p.~99]{Andrews}, for each $k \geq 1$, the element $2$ generates $(\Z/p^k \Z)^\times$, so the Frobenius automorphism acts transitively on the roots of $\Phi_{p^k}$ in $\Fbar_2$, so $\Phi_{p^k}$ is irreducible over $\FF_2$.

We have
\begin{equation} \label{eq:f2Phi}
f_2 = X^n + X^{n-1} + \ldots + 1 = \Phi_p \Phi_{p^2} \cdot \ldots \cdot \Phi_{p^r} \in \FF_2[X]
\end{equation}
and by the above argument each of the cyclotomic polynomials on the right-hand side in \eqref{eq:f2Phi} is irreducible. 
If $f$ is reducible over $\Q$, then $f=gh$ for  some monic $g,h\in \Z[X]$ of positive degrees. Then  $f_2=g_2h_2$, so $\Phi_{p^r}$ divides one of the factors over $\FF_2$, say $g_2$, and then 
\begin{equation} \label{eq:divdeg}
  \deg h = \deg f - \deg g \leq \deg f - \deg \Phi_{p^r} = (p^r-1) - (p^r - p^{r-1}) = p^{r-1} - 1 < n/p < n/1470.
\end{equation}

On the other hand, we will use \cite[Theorem~7]{BKK} to prove that if $n$ is sufficiently large, then
\begin{equation}
\label{E:probability of divisor of degree <= n/1470}
\textup{$\mb{P}(\textup{$f$ has a divisor of degree $\leq n/1470$})$ is $O(n^{-c})$.}
\end{equation}
For $\gamma = 1/2$, $P = 3 \cdot 5 \cdot 7 \cdot 11 = 1155$, and $s = 735$, we verify numerically that for all integers $Q,R,\ell$ with $QR=P$ and $Q>1$, 
\[
     \sum_{k \in \Z/Q\Z} \abs{\hat{\mu}(k/Q+\ell/R)}^s \le 0.9999 \, Q^{1-\gamma};
\]
thus if $n$ is sufficiently large and $\mu_j \coloneqq \mu$ for $j=0,\ldots,n-1$, then the conditions of \cite[Theorem~7]{BKK} are satisfied, so \eqref{E:probability of divisor of degree <= n/1470} holds. 
By \eqref{eq:divdeg} and~\eqref{E:probability of divisor of degree <= n/1470}, with probability at least $1-O(n^{-c})$, we obtain a contradiction showing that $f$ is irreducible.
\end{proof}

\section{Approximate equidistribution and divisors modulo \texorpdfstring{$3$}{3}}
\label{sec:mod3}

The main difference between the proof of Theorem~\ref{thm:main} for large primes in the previous section \S \ref{sec:largep} and the one for $7 \le p \le 1470$ in \S \ref{sec:mainproof} is that the black box \cite[Theorem~7]{BKK} in the former combines information on the factorization of $f_3$, $f_5$, $f_7$ and $f_{11}$ to rule out \emph{all} positive integers $\leq n/1470$ as possible degrees of divisors of $f$; in contrast, since the factorization of $f_2$ shows that already many such values cannot occur as degrees of a divisor, in \S \ref{sec:mainproof} we require information on the factorization of $f_3$ only. 

The aim of this section is to bound the probability that the degree~$n$ polynomial $f_3$ has a divisor of fixed degree $k$.  For this, we require some further results of \cite{BKK}. Define
\begin{equation*}
\Delta_p(n; m) \coloneqq \sum_{\substack{D \in \FF_p[X], \, X \nmid D \\ \deg{D} \leq m}} \max_{C \in \FF_p[X]} \abs{ \mb{P}(f_p \equiv C \bmod{D}) - p^{-\deg{D}} }.
\end{equation*}
The quantity $\Delta_p$ measures the extent to which $f_p$ fails to be equidistributed modulo $D \in \FF_p[X]$ on average for all $D$ not divisible by $X$ of degree at most $m$. We study it for $p = 3$.
Set $\theta^* \coloneqq \frac{\log 2}{2\log 3} \approx 0.315$. 

\begin{lemma} \label{lem:Delta}
Fix a positive real number $\theta < \theta^*$.  Then, as $n \to \infty$,  
\begin{equation*}
    \Delta_3(n; \theta n) \ll e^{-n^{1/10}}.
\end{equation*}
\end{lemma}

\begin{proof}
Applying \cite[Proposition~2.3]{BKK} with $\mc{P} = \{3\}$, we have 
\begin{equation*}
\Delta_3(n; \gamma n/s + n^{0.88}) \ll e^{-n^{1/10}}
\end{equation*}
for any $\gamma \geq 1/2$, positive integer $s$, and sufficiently large $n$ satisfying
\begin{equation} \label{eq:fouriersum}
\abs{\hat{\mu}(0)}^s + \abs{\hat{\mu}(1/3)}^s + \abs{\hat{\mu}(2/3)}^s \leq (1-n^{-1/10}) \cdot 3^{1-\gamma}.
\end{equation}
Using \eqref{eq:FT}, we find that \eqref{eq:fouriersum} holds for any $\gamma < \gamma(s) \coloneqq 1-\frac{\log(1+2^{1-s})}{\log 3}$ and $n$ sufficiently large (depending on $\gamma$). The smallest $s$ such that $\gamma(s) > 1/2$ is $s = 2$, in which case we find $\gamma(s)/s = \theta^*$.
\end{proof}

Let $\tau(f_p)$ be the number of (not necessarily irreducible) monic divisors in $\FF_p[X]$ of the polynomial $f_p$. Denote by $f_p^{\mc{S}(m)}$ the $m$-smooth part of $f_p$, defined as the product of all monic irreducible divisors in $\FF_p[X]$ of $f_p$ (with multiplicity) of degree $\leq m$.

\begin{definition}
\label{def:b_smooth}
Fix $\epsilon \in (0,1/2)$ and a positive real number $\theta<\theta^*$. Let $n \in \Z_{\ge 1}$ and $k \in \R_{\geq 1}$.
For a random Littlewood polynomial $f$ of degree $n$, define the \emph{bounded smoothness event} $\mc{E}_{k, \theta, \epsilon,n}$ as the event in which
\begin{equation*}
\deg(f_3^{\mc{S}(m)}) \leq \epsilon m \log m \quad \text{and} \quad \tau(f_3^{\mc{S}(m)}) \leq m^{(1+\epsilon)\log 2}
\end{equation*}
holds for all integers $m \in [k, 2\theta n/\log n]$.
\end{definition}

\begin{lemma}
\label{lem:crux}
Fix $\epsilon\in (0,1/2)$ and $\theta<\theta^*$.  Then there exist $c,C,C'>0$ such that for all sufficiently large $n$ and all positive integers $k < \theta n$, we have 
    \begin{itemize}
    \item[(i)] $\mb{P}(\mc{E}_{k, \theta, \epsilon,n})  > 1- C k^{-c}$, and 
    \item[(ii)]  $ \mb{P}(\mc{E}_{k^{1/4},\theta, \epsilon,n} \text{ holds and } f_3 \text{ has a divisor of degree } k) < C' \dfrac{(\log n)^2}{k^{1-\log 2-2\epsilon}}$.
    \end{itemize}
\end{lemma}

\begin{proof}\hfill
\begin{itemize}
\item[(i)]
Apply Lemma~\ref{lem:Delta} and \cite[Lemma~9.3]{BKK}, where $m_0 = k$, $\mc{P} = \{3\}$, and $\theta$ and $\epsilon$ are as here. 
\item[(ii)]
We will use \cite[Lemma~9.4]{BKK}. In \cite{BKK} there is an extra parameter, $\lambda$, which will not be important for us; we set $\lambda = 1-\epsilon$. The condition $\Delta_3(n;\theta n+n^\lambda)\leq n^{-7}$ of \cite[Lemma~9.4]{BKK} follows by using Lemma~\ref{lem:Delta} with $\theta_{\textrm{Lemma~\ref{lem:Delta}}}=\frac12 (\theta+\theta^*)$, if $n$ is sufficiently large (depending on $\theta$ and $\epsilon$). Further parameters of \cite[Lemma~9.4]{BKK} are $\delta = 1$ and again $\mc{P} = \{3\}$, and $\theta$, $\epsilon$, and $k$ as here. The $\mc{E}_{k, \lambda, \epsilon, \theta}$ of \cite{BKK} is a larger event than our $\mc{E}_{k^{1/4}, \theta, \epsilon, n}$, since our $k^{1/4}$ is a lower bound for the $m_0$ in \cite[Lemma~9.4]{BKK}. 
Thus \cite[Lemma~9.4]{BKK} implies the desired bound, but with $k^{(1-\log 2 -\epsilon)\lambda}$ in the denominator. This is stronger than needed, since $(1-\log 2-\epsilon)\lambda > 1-\log 2-2\epsilon$.\qedhere
\end{itemize}
\end{proof}

\section{Proof of Theorem~\ref{thm:main} for $7 \le p \le 1470$}
\label{sec:mainproof}

Since $2$ does not generate $(\Z/7^2\Z)^\times$, we have $p\geq 11$.
We fix $p$ and assume throughout that $n$ is sufficiently large.
For $1 \le k < n$, let $E_k$ be the event that $f$ has a degree~$k$ factor.
For $I \subset \R$, let $E_I = \Union_{k \in I} E_k$.
By the argument of \S \ref{sec:largep}, if $f$ is reducible, it has a factor $h$ such that $h_2$ is a subproduct of $\Phi_p \Phi_{p^2} \cdot \ldots \cdot \Phi_{p^{r-1}}$.
For $j \le r-2$, let $D_j$ be the set of integers $k>n^{1/10}$ such that $k$ is the degree of such a subproduct whose largest factor is $\Phi_{p^{j+1}}$.
Then $\# D_j \leq 2^j$, and $D_j \subset [p^j,p^{j+1})$.
Let $s$ be the largest positive integer such that $p^s \le n^{1/10}$. 
Then the event that $f$ is reducible is contained in $E_{[1,n^{1/10}]}  \union \Union_{j=s}^{r-2} E_{D_j}$.

Let $\theta \coloneqq 1/p < \frac{\log 2}{2\log 3} = \theta^*$.
Let $\epsilon = 0.001$.
We have $n^{1/40} < \theta n$ for large $n$;
let $\mc{E} = \mc{E}_{n^{1/40}, \theta, \epsilon, n}$ (see Definition~\ref{def:b_smooth}). Let $\mc{E}^{\mr{comp}}$ be the complementary event.
Then
\begin{equation} \label{E:P(f is reducible)}
\mb{P}(\textup{$f$ is reducible}) \le \mb{P}(E_{[1,n^{1/10}]}) + \mb{P}(\mc{E}^{\mr{comp}}) + \sum_{j=s}^{r-2} \sum_{k \in D_j} \mb{P}(\mc{E} \intersect E_k).
\end{equation}

By \cite[Proposition~2.1]{BKK}, $\mb{P}(E_{[1,n^{1/10}]}) \ll n^{-7/20}$. 
By Lemma~\ref{lem:crux}(i), $\mb{P}(\mc{E}^{\mr{comp}}) \ll n^{-c'}$ for some $c'>0$.
For $n^{1/10} < k < \theta n$, we have $n^{1/40} < k^{1/4} < \theta n$; then $\mc{E} \subset \mc{E}_{k^{1/4}, \theta, \epsilon,n}$, and $E_k$ is contained in the event that $f_3$ has a degree~$k$ divisor, so 
\begin{align*}
   \mb{P}(\mc{E} \intersect E_k) &\le \mb{P}(\textup{$\mc{E}_{k^{1/4}, \theta, \epsilon,n}$ holds and $f_3$ has a degree~$k$ divisor}) \\
   &\ll \dfrac{(\log n)^2}{k^{0.3}} \qquad\textup{(by Lemma~\ref{lem:crux}(ii), since $1-\log 2 - 2\epsilon \ge 0.3$);} \\
   \sum_{j=s}^{r-2} \sum_{k \in D_j} \mb{P}(\mc{E} \intersect E_k)
        &\ll \sum_{j=s}^{r-2} \#D_j \dfrac{(\log n)^2}{p^{0.3 j}} 
        \leq (\log n)^2 \sum_{j=s}^{r-2} \left( \dfrac{2}{p^{0.3}} \right)^j \ll (\log n)^2 p^{-0.01 s}
        \ll n^{-c''}
\end{align*}
for any positive $c'' < 0.001$, since $2/p^{0.3} \le p^{-0.01}$ for $p \ge 11$, and $p^s \gg p^{s+1} \ge n^{1/10}$ since $p$ was fixed.
We have now bounded all three terms on the right of \eqref{E:P(f is reducible)}, so $\mb{P}(\textup{$f$ is reducible}) \ll n^{-c}$ for some $c>0$.

\bibliographystyle{amsplain}

\end{document}